\documentclass[12pt,a4paper]{article}
\usepackage{graphicx}
\usepackage{amsmath}
\usepackage{amssymb}
\usepackage{amsthm}
\usepackage{geometry}
\usepackage{fancyhdr}
\usepackage{color} 
\usepackage{overpic}
\usepackage{mathabx}

\usepackage[all]{xy}
\newcommand{\al}{\alpha}

\newcommand{\be}{\beta}
\newcommand{\ga}{\gamma}

\newcommand{\ra}{\rightarrow}

\newcommand{\bpf}{\begin{proof}}
\newcommand{\epf}{\end{proof}}
\newcommand{\bthm}{\begin{thm}}
\newcommand{\ethm}{\end{thm}}
\newcommand{\bprop}{\begin{prop}}
\newcommand{\eprop}{\end{prop}}
\newcommand{\bcor}{\begin{cor}}
\newcommand{\ecor}{\end{cor}}
\newcommand{\blem}{\begin{lem}}
\newcommand{\elem}{\end{lem}}
\newcommand{\bdefn}{\begin{defn}}
\newcommand{\edefn}{\end{defn}}
\newcommand{\bexmp}{\begin{exmp}}
\newcommand{\eexmp}{\end{exmp}}
\newcommand{\brem}{\begin{rem}}
\newcommand{\erem}{\end{rem}}

\newcommand{\bdia}{\begin{displaymath}\xymatrix}
\newcommand{\edia}{\end{displaymath}}
\newcommand{\beq}{\begin{equation*}\begin{aligned}}
\newcommand{\eeq}{\end{aligned}\end{equation*}}

\newcommand{\intg}{\mathbb{Z}}
\newcommand{\real}{\mathbb{R}}
\newcommand{\comp}{\mathbb{C}}

\newtheorem{thm}{\textbf {Theorem}}[section]
\newtheorem{cor}[thm]{\textbf{Corollary}}
\newtheorem{prop}[thm]{\textbf{Proposition}}
\newtheorem{lem}[thm]{\textbf{Lemma}}
\newtheorem{conj}[thm]{Conjecture}

\theoremstyle{definition}
\newtheorem{defn}[thm]{\textbf{Definition}}

\newtheorem{exmp}[thm]{Example}

\theoremstyle{remark}
\newtheorem{rem}[thm]{Remark}

\newcommand{\Ga}{\Gamma}

\newcommand{\shi}{\underline{\rm SHI}}

\title{Instanton and the depth of taut foliations}

\author{Zhenkun Li}
\date{}

\begin{document}
\bibliographystyle{plain}

\maketitle
\abstract{Sutured instanton Floer homology was introduced by Kronheimer and Mrowka in \cite{kronheimer2010knots}. In this paper, we prove that for a taut balanced sutured manifold with vanishing second homology, the dimension of the sutured instanton Floer homology provides a bound on the minimal depth of all possible taut foliations on that balanced sutured manifold. The same argument can be adapted to the monopole and also the Heegaard Floer settings, which gives a partial answer to one of Juhasz's conjectures in \cite{juhasz2008floer}. Also, using the nature of instanton Floer homology, on knot complements, we can construct a taut foliation with bounded depth, given some information on the representation varieties of the knot fundamental groups. This indicates a mystery relation between the representation varieties and some small depth taut foliations on knot complements and gives a partial answer to one of Kronheimer and Mrowka's conjecture in \cite{kronheimer2010knots}.


\section{Introduction}
\subsection{Statement of results}
Sutured manifolds and sutured manifold hierarchies were introduced by Gabai in 1983 in \cite{gabai1983foliations} and subsequent papers. They are powerful tools in the study of $3$-dimensional topology. A sutured manifold is a compact oriented $3$-manifold with boundary, together with an oriented closed $1$-submanifold $\ga$ on $\partial{M}$, which is called the suture. If $S\subset M$ is a properly embedded surface inside $M$, which satisfies some mild conditions, then we can perform a decomposition of $(M,\ga)$ along $S$ and obtain a new sutured manifold $(M',\ga')$. We call this process a sutured manifold decomposition and write
$$(M,\ga)\stackrel{S}{\leadsto}(M',\ga').$$

A balanced sutured manifold $(M,\ga)$ is a sutured manifold with some further restrictions on $M$ and $\ga$. It was introduced by Juh\'asz in \cite{juhasz2006holomorphic} to accommodate the construction of Heegaard Floer homology on them. Later, when Kronheimer and Mrowka introduced sutured monopole and instanton Floer homologies, they also used the settings of balanced sutured manifolds. So, in this paper, we will only work with balanced sutured manifolds, though it should be understood that Gabai's results were initially proved for general sutured manifolds.

A celebrated theorem proved by Gabai is the following.
\bthm[Gabai \cite{gabai1983foliations}]\label{thm: hierarchy}
Suppose $(M,\ga)$ is a taut balanced sutured manifold, then there exists a finite sequence of sutured manifold decompositions
\begin{equation}\label{eq: hierary}
(M,\ga)\stackrel{S_1}{\leadsto}(M_1,\ga_1)\stackrel{S_2}{\leadsto}...\stackrel{S_n}{\leadsto}(M_n,\ga_n),
\end{equation}
where $(M_i,\ga_i)$ is taut for all $i$ and $(M_n,\ga_n)$ is a product sutured manifold, meaning that there is an oriented surface $F$ with non-trivial boundary so that 
$$(M_n,\ga_n)=([-1,1]\times F,\{0\}\times \partial{F}).$$
\ethm 

One original motivation for Gabai to establish Theorem \ref{thm: hierarchy} is to construct taut foliations on $3$-manifolds. In particular, he proved the following theorem.

\bthm[Gabai \cite{gabai1983foliations}]\label{thm: existence of finite depth taut foliations}
Suppose $(M,\ga)$ is a taut balanced sutured manifold, then $(M,\ga)$ admits a finite depth taut foliation.
\ethm

However, Gabai only proved the existence of a taut foliation with finite depth, yet he didn't offer any bounds on how small the depth could be. In \cite{juhasz2010polytope}, Juh\'asz  made the following conjecture.
\begin{conj}[Juh\'asz \cite{juhasz2010polytope}]\label{conj: juhasz}
Suppose $(M,\ga)$ is a taut balanced sutured manifold with $H_2(M)=0$, and
$${\rm rk}_{\intg_2}(SFH(M,\ga))<2^{k+1},$$
then $(M,\ga)$ admits a taut foliation of depth at most $2k$.
\end{conj}

Here, $SFH(M,\ga)$ is the sutured (Heegaard) Floer homology of $(M,\ga)$, introduced by Juhasz in \cite{juhasz2006holomorphic}. It is a finite-dimensional vector space over the field $\intg_2$ and is a topological invariant associated to $(M,\ga)$. Following this line, Kronheimer and Mrowka further made the following conjecture.
\begin{conj}[Kronheimer and Mrowka \cite{kronheimer2010knots}]\label{conj: KM}
    Suppose $K\subset S^3$ is a knot, and consider the irreducible homomorphisms
    $$\rho:\pi_1(S^3(K)\ra SU(2))$$
    which map a chosen meridian $m$ to the element $\mathbf{i}\in SU(2)$. Suppose that these homomorphisms are non-degenerate and that the number of conjugacy classes of such homomorphisms is less than $2^{k+1}$. Then, the knot complement $S^3(K)$ admits a taut foliation of depth at most $2k$, transverse to the boundary $\partial{S^3(K)}$.
\end{conj}

In this paper, we prove the following result, constructing a taut foliation whose depth is bounded in terms of the dimension of sutured instanton Floer homology. The sutured instanton Floer homology, denoted by $SHI$, is another type of Floer theory associated to $(M,\ga)$, which was introduced by Kronheimer and Mrowka in \cite{kronheimer2010knots}. It is a finite-dimensional vector space over $\mathbb{C}$. Though the bound on the depth of taut foliations is not as sharp as the ones in Conjecture \ref{conj: juhasz} and Conjecture \ref{conj: KM}, up to the author's knowledge, it is the first of this kind.

\bthm\label{thm: small depth taut foliations}
Suppose $(M,\ga)$ is a taut balanced sutured manifold with $H_2(M)=0$, and 
$${\rm dim}_{\mathbb{C}}SHI(M,\ga)<2^{k+1}.$$
Then, $(M,\ga)$ admits a taut foliation of depth at most $2^{k+6}$.
\ethm

\bcor\label{cor: small depth taut foliations}
Conjecture \ref{conj: juhasz} and \ref{conj: KM} hold if we replace the depth $2k$ of the taut foliation in the statement of the conjecture by $2^{k+6}$.
\ecor

\brem
It can be reformulated that the minimal depth of all taut foliations is bounded by a multiple of the dimension of the sutured instanton Floer homology. However, the original statement in Theorem \ref{thm: small depth taut foliations} is more convenient for carrying out the proof.
\erem

\subsection{Strategy of the proof}
In this subsection, we roughly explain the technical difficulty of attacking Conjecture \ref{conj: juhasz} and the idea to obtain the main result of the current paper. Though the original conjecture was stated in the Heegaard Floer setup, we present everything in the instanton settings, in consistence with the main proof below. 

Suppose we have a taut balanced sutured manifold $(M,\ga)$ with $H_2(M)=0$ and whose instanton Floer homology has dimension smaller than $2^{k+1}$. Then, one can construct a sutured manifold hierarchy in $k$ stage so that in each stage, we decompose the sutured manifold twice, and that in each stage, the dimension of the Floer homology of the resulting balanced sutured manifold after the two decompositions is at most a half of the previous one. See Juh\'asz \cite{juhasz2010polytope} or Ghosh and Li \cite{li2019decomposition}.

Suppose, for simplicity, that all balanced sutured manifolds involved are horizontally prime. Then, in each stage, the first of the two decompositions is to make the balanced sutured manifold free of essential product annuli and essential product disks so that one could apply Proposition \ref{prop: decomposition drops the dimension by half} to perform the second decomposition. Thus, after $k$ stages and altogether $2k$ decompositions, we obtain a product sutured manifold which admits a taut foliation of depth $0$. By Theorem \ref{thm: well groomed decomposition increase depth by 1}, if the each one of the $2k$ decompositions is well-groomed, then we can glue the taut foliation on the product sutured manifold through the hierarchy and obtain a depth-$2k$ taut foliation on $(M,\ga)$, thus proving Conjecture \ref{conj: juhasz}.

In each stage, the second of the two decompositions is coming from Proposition \ref{prop: decomposition drops the dimension by half} and is well-groomed already. However, the first of the two decompositions is to decompose along a maximal disjoint union of essential product annuli, which is not necessarily well-groomed. So a naive attempt to attack Conjecture \ref{conj: juhasz} fails at this point. In this paper, instead of showing that the first of the two decompositions can be made well-groomed, we show that it can be decomposed into a sequence of decompositions where each decomposing surface is either well-groomed or having its boundary contained in an annular neighborhood of the suture. Furthermore, each of such decompositions would only increase the depth of the taut foliation we glue back by at most one, and the total number of such decompositions can be bounded above by the dimension of the Floer homology of the sutured manifold we start with.

To obtain this bound, recall that in each stage, the first of the two decompositions is to decompose along a maximal disjoint union of essential product annuli. The decomposition along one essential product annulus can be decomposed into a sequence of decompositions where each decomposing surface is either well-groomed or having its boundary contained in an annular neighborhood of the suture. Furthermore, the total number of such decompositions is a fixed constant. So we only need to bound the number of product annuli that are involved. This is done through Proposition \ref{prop: bounding dimension by genus of the boundary}: if the balanced sutured manifold is free of essential product annuli and product disks (and we have assumed it to be horizontally prime), then the dimension of Floer homology bounds the genus of the boundary of the $3$-manifold and hence the number of components of sutures, which further gives a bound on the number of product annuli we need to decompose along, thus completing the proof.

{\bf Acknowledgement}. This material is based upon work supported by the National Science Foundation under Grant No. 1808794. The author would like to thank his advisor Tomasz S. Mrowka for his enormous help and valuable suggestions. The author would also like to thank the referee for helpful comments.

\section{Preliminaries}
In this paper, all the notations will be kept the same as in Ghosh and Li \cite{li2019decomposition}. So if a term has already been defined in that paper, we will not define it again. We have the following new definitions.

\bdefn[Gabai \cite{gabai1983foliations}]\label{defn: well groomed surface}
Suppose $(M,\ga)$ is a balanced sutured manifold. A surface $S$ is called {\it well-groomed} if $\partial S$ is essential in $H_1(\partial M)$, and the following is true.

(1) For each component $A$ of the annular neighborhood $A(\ga)$, $S\cap A$ consists of either a collection of parallel and coherently oriented non-separating simple arcs or a collection of parallel simple closed curves each oriented in the same way as $\ga\cap A$.

(2) For each component $V$ of $R(\ga)$, $S\cap V$ consists of either a collection of parallel and coherently oriented non-separating simple arcs or a collection of parallel and coherently oriented non-separating simple closed curves.
\edefn

\bdefn[Gabai \cite{gabai1983foliations}]\label{defn: taut foliations}
A transversely oriented co-dimension-one foliation $\mathfrak{F}$ on a balanced sutured manifold $(M,\ga)$ is called {\it taut} if $\mathfrak{F}$ is transverse to $A(\ga)$, tangent to $R(\ga)$ with normal direction pointing inward on $R_-(\ga)$ and outward on $R_+(\ga)$, $\mathfrak{F}|_{A(\ga)}$ has no Reeb components, and there exists a (not necessarily connected and not necessarily closed) curve $c$ so that each leaf of $\mathfrak{F}$ intersects $c$ transversely and non-trivially.
\edefn

\bdefn[Gabai \cite{gabai1983foliations}]\label{defn: depth of foliations}
Let $M$ be a compact oriented $3$-manifold, and $\mathfrak{F}$ a co-dimension-one foliation. We say a leaf $L$ of $\mathfrak{F}$ has {\it depth} $0$ if $L$ is compact. Suppose we have defined the depth for $j\leq k$, then say a leaf $L$ of $\mathfrak{F}$ has {\it depth} $k+1$ if $\widebar{L}\backslash L$ is a union of leaves of depth at most $k$ and contains at least one leaf of depth $k$. The foliation $\mathfrak{F}$ is called {\it depth} $k$ if all its leaves have depth at most $k$, and it admits at least one leaf of depth $k$. If such a $k$ does not exist, then we say $\mathfrak{F}$ has infinite depth.
\edefn

The following two Propositions from Ghosh and Li \cite{li2019decomposition} are key ingredients of the proof of Theorem \ref{thm: small depth taut foliations}.

\bprop\label{prop: decomposition drops the dimension by half}
Suppose $(M,\ga)$ is a connected balanced sutured manifold so that $H_2(M)=0$ and it is horizontally prime, taut, and free of essential product annuli and product disks. If $(M,\ga)$ is not a product sutured manifold, then there exists a well-groomed surface $S\subset M$ and a sutured manifold decomposition
$$(M,\ga)\stackrel{S}{\leadsto}(M',\ga')$$
so that $(M',\ga')$ is taut, and
$${\rm dim}_{\mathbb{C}}SHI(M',\ga')\leq\frac{1}{2}{\rm dim}_{\mathbb{C}}SHI(M,\ga).$$
\eprop

\bpf
This is essentially \cite[Proposition 5.8]{li2019decomposition}. In the original statement of the proposition, the homology class $\al\in H_2(M,\partial M)$ can be chosen freely, so we can choose a well-groomed class as guaranteed by \cite[Lemma 2.8]{juhasz2010polytope}. Then, Proposition \ref{prop: decomposition drops the dimension by half} follows.
\epf

\bprop\label{prop: bounding dimension by genus of the boundary}
Suppose $(M,\ga)$ is a connected balanced sutured manifold so that $H_2(M)=0$ and it is horizontally prime, taut, and free of essential product annuli and product disks. Then, 
$${\rm dim}_{\mathbb{C}}SHI(M,\ga)\geq g(\partial{M})+1.$$
\eprop
\bpf
By Corollary 5.12 in Ghosh and Li \cite{li2019decomposition}, we know that
$${\rm dim}_{\mathbb{C}}SHI(M,\ga)\geq {\rm dim}_{\real}H^2(M,\partial M;\real)+1.$$
Since $H_2(M)=0$ and $M$ is connected, we know that
$${\rm dim}_{\real}H^2(M,\partial M;\real)=g(\partial{M}).$$
\epf

To construct a taut foliation with controlled depth, we also need the following theorem from Gabai \cite{gabai1983foliations}.

\bthm\label{thm: well groomed decomposition increase depth by 1}
Suppose $(M,\ga)$ is a balanced sutured manifold. Suppose $S\subset M$ is a well-groomed surface so that we have a decomposition
$$(M,\ga)\stackrel{S}{\leadsto}(M',\ga').$$
Suppose further that $(M',\ga')$ admits a taut foliation of depth $k$, then $(M,\ga)$ admits a taut foliation of depth at most $k+1$.
\ethm

The following Proposition is essentially \cite[Proposition 8.10]{juhasz2008floer} and is adapted to the instanton theory by Kronheimer and Mrowka.

\begin{prop}[Kronheimer and Mrowka {\cite[Proposition 6.7]{kronheimer2010knots}}]\label{prop: decompose along product annuli do not increase SHI}
Suppose $(M,\gamma)$ is a taut balanced sutured manifold so that $H_2(M)=0$. Suppose $A\subset (M,\gamma)$ is a non-trivial product annulus. Suppose further $(M',\gamma')$ is obtained from $(M,\gamma)$ by decomposing along $A$. Then we have
$${\rm dim}_{\comp}SHI(M',\ga')\leq {\rm dim}_{\comp}SHI(M,\ga).$$
\end{prop}

In the next section, the following definition is also useful. 
\bdefn\label{defn: quasi-horizontal surfaces}
Suppose $(M,\ga)$ is a balanced sutured manifold. A properly embedded surface $S$ is called {\it quasi-horizontal} if $\partial S\cap R(\ga)=\emptyset$ and for any component $A$ of $A(\gamma)$, $S\cap A$ is either empty or parallel simple closed curves each parallel to $\ga$ as oriented curves. 
\edefn

Suppose $(M,\ga)$ is a balanced sutured manifold and $S\subset (M,\ga)$ is a quasi-horizontal surface. Suppose further that $(M',\ga')$ is obtained from $(M,\ga)$ by decomposing along $S$, then we know that $R_+(\ga')$ has a component $V_+$ and $R_-(\ga')$ has a component $V_-$, so that $V_+$ and $V_-$ are both diffeomorphic to $S$ and $(M,\gamma)$ is obtained from $(M',\ga')$ by identifying $V_+$ with $V_-$ via the diffeomorphisms to $S$. The following lemma is straightforward. 
\blem\label{lem: glue taut foliation through quasi-horizontal surfaces}
Suppose $(M',\gamma')$ is a taut balanced sutured manifold admitting a finite depth taut foliation $\mathfrak{F}'$. Suppose $V_+$ is a component of $R_+(\ga')$ and $V_-$ is a component of $R_-(\ga')$. Suppose further that there is an orientation preserving diffeomorphism $f:V_+\ra V_-$. We can glue $V_+$ to $V_-$ via $f$ to obtain a new sutured manifold $(M,\ga)$. Then there exists a taut foliation $\mathfrak{F}$ on $(M,\gamma)$ whose depth is at most that of $\mathfrak{F}'$.
\elem

\section{Constructing taut foliations of bounded depth}
\blem\label{lem: curves on surface}
Suppose $\Sigma$ is a closed connected oriented surface of genus $g\geq 1$. Then there are at most $(3g-2)$ many connected simple closed curves on $\Sigma$ so that they are each non-separating, pairwise disjoint, and pairwise non-parallel.
\elem
\bpf
Suppose $\ga$ is a collection of connected simple closed curves on $\Sigma$ so that they are each non-separating, pairwise disjoint, and pairwise non-parallel. Write $|\ga|$ to be the number of components of $\ga$. We want to show that $|\ga|\leq 3g-2$.

First assume $g>1$. Write $S=\Sigma\backslash N(\ga)$, where $N(\ga)$ is an annular neighborhood of $\ga$. By assumption, each component of $S$ has negative Euler Characteristics. Since $\chi(S)=\chi(\Sigma)=2-2g,$ we know that $|S|$, i.e., the number of components of $S$ is at most $2g-2$. Also, we know that
$$|\partial{S}|-2\cdot |S|\leq -\chi(S)=2g-2,$$
so
$$|\ga|=\frac{1}{2}|\partial S|\leq 3g-3.$$

When $g=1$, clearly $|\ga|\leq 1$ and we are done.
\epf


Now we are ready to prove the main result of the paper.
\bpf[Proof of Theorem \ref{thm: small depth taut foliations}]
We prove the theorem by the induction on $k$. When $k=0$, by \cite[Theorem 1.2]{li2019decomposition}, we know that $(M,\ga)$ is a product sutured manifold, and hence it admits a taut foliation of depth $0$. Suppose the theorem holds for $k<k_0$. Now we argue for the case $k=k_0$. (We will keep writing $k$ instead of $k_0$)

{\bf Case 1}. The balanced sutured manifold $(M,\ga)$ is not horizontally prime. Then, we can find a non-boundary-parallel horizontal surface $S$ and perform a decomposition:
$$(M,\ga)\stackrel{S}{\leadsto}(M_{1},\ga_{1}),$$
with $(M_{1},\ga_{1})$ being a disjoint union
$$(M_{1},\ga_{1})= (M_{2},\ga_{2})\sqcup(M_{3},\ga_{3}).$$
Since $S$ is not boundary parallel, we conclude that
$${\rm dim}_{\mathbb{C}}SHI(M_{2},\ga_2)\leq \frac{1}{2}{\rm dim}_{\mathbb{C}}SHI(M,\ga)$$
and
$${\rm dim}_{\mathbb{C}}SHI(M_3,\ga_3)\leq \frac{1}{2}{\rm dim}_{\mathbb{C}}SHI(M,\ga).$$
Hence, by the inductive hypothesis, $(M_2,\ga_2)$ and $(M_3,\ga_3)$ both admits a taut foliation of depth $2^{k+5}$. Since $S$ is quasi-horizontal, we are done by Lemma \ref{lem: glue taut foliation through quasi-horizontal surfaces}.

{\bf Case 2}. The balanced sutured manifold $(M,\ga)$ is horizontally prime. According to the proof of \cite[Proposition 7.6]{juhasz2010polytope} (also c.f. \cite[Proposition 2.16]{juhasz2010polytope} and \cite[Lemma 4.2]{scharlemann2006three}), we know that there is a disjoint union $A$ of non-trivial product annuli, and a sutured manifold decomposition
$$(M,\ga)\stackrel{A}{\leadsto}(M',\ga'),$$
so that $(M',\ga')$ is taut, horizontally prime, reduced, and with $H_2(M')=0$. Let $(M_1,\ga_1)$ be the union of components of $(M',\ga')$ that are not product sutured manifolds. There is a minimal union of product annuli $A'\subset A$ so that the decomposition along $A'$ results in
$$(M,\ga)\stackrel{A'}{\leadsto}(M_{1},\ga_{1})\sqcup(M_2,\ga_2),$$
where $(M_2,\ga_2)$ is a product sutured manifold. Write
$$(M_3,\ga_3)=(M_1,\ga_1)\cup (M_2,\ga_2),$$
we know that $(M_1,\ga_1)$ is also taut, horizontally prime, reduced, and with $H_2(M_1)=0$.

Suppose the components of $A'$ are
\begin{equation}\label{eq: definition of n}
A'=A_1\cup...\cup A_n.
\end{equation}

{\bf Claim 1}. We have
$$n<6\cdot 2^{k+1}.$$

To prove this claim, by Proposition \ref{prop: decompose along product annuli do not increase SHI} it suffices to show that
$$n<6\cdot {\rm dim}_{\mathbb{C}}(\shi(M_1,\ga_1))$$

First assume that $(M_1,\ga_1)$ is connected. By Proposition \ref{prop: bounding dimension by genus of the boundary}, we know that
$$g(\partial M_1)< {\rm dim}_{\mathbb{C}}(\shi(M_1,\ga_1)).$$
Assume that $|\ga_1|>6\cdot {\rm dim}_{\mathbb{C}}(\shi(M_1,\ga_1))-4$, then by Lemma \ref{lem: curves on surface} and the pigeon-hole theorem, we know that there are three components of $\ga_1$ that are parallel to each other. Hence, there is a non-trivial product annulus separating the three parallel sutures from the rest, which contradicts the fact that $(M_1,\ga_1)$ is reduced. Thus, we conclude that
$$n<|\ga_1|\leq 6\cdot {\rm dim}_{\mathbb{C}}(\shi(M_1,\ga_1))-4<6\cdot {\rm dim}_{\mathbb{C}}(\shi(M_1,\ga_1)).$$

In general, if $(M_1,\ga_1)$ is disconnected, then we can proceed by a second induction. Suppose we have proved that
$$n<6\cdot {\rm dim}_{\mathbb{C}}(\shi(M_1,\ga_1))$$
if $(M_1,\ga_1)$ has $m$ components. If $(M_1,\ga_1)$ has $m+1$ components, then assume
$$(M_1,\gamma_1)=(M_{1,1},\ga_{1,1})\cup(M_{1,2},\ga_{1,2}),$$
where $(M_{1,1},\ga_{1,1})$ has $m$ component and $(M_{1,2},\ga_{1,2})$ is connected. Then by the inductive hypothesis, we know that
\beq
|\ga_1|<&6\cdot {\rm dim}_{\mathbb{C}}(\shi(M_{1,1},\ga_{1,1}))+6\cdot {\rm dim}_{\mathbb{C}}(\shi(M_{1,2},\ga_{1,2}))\\
<&6\cdot {\rm dim}_{\mathbb{C}}(\shi(M_1,\ga_1)).
\eeq
The last equality holds, since by \cite[Proposition 6.5]{kronheimer2010knots} we know that
$$\shi(M_1,\ga_1)\cong\shi(M_{1,1},\ga_{1,1})\otimes\shi(M_{1,2},\ga_{1,2}),$$
and by construction, both $(M_{1,1},\ga_{1,1})$ and $(M_{1,2},\ga_{1,2})$ are not product and hence have Floer homology of dimension at least two.

{\bf Claim 2}. There is a sequence of sutured manifold decompositions
$$(M,\ga)\stackrel{S_1}{\leadsto}(N_{1},\delta_{1})...\stackrel{S_n}{\leadsto}(N_{n},\delta_{n}),$$
where $n$ is defined as in Formula (\ref{eq: definition of n}), so that the following is true.

(1) Each $S_i$ is either well-groomed or quasi-horizontal.

(2) Each $(N_i,\delta_i)$ is taut.

(3) Suppose the components of $\ga_3$ are
$$\ga_3=\theta_1\cup...\cup\theta_m.$$
Then, for each $i=1,...,m$, there is a compact connected oriented surface-with-boundary $F_i$ satisfying the following properties.

(a) For $i=1,...,m$, there is an orientation reversing embedding
$$f_i:\theta_i\hookrightarrow \partial{F}_i.$$

(b) Write 
$$F=F_1\cup...\cup F_m ~{\rm and}~ f=f_1\cup...\cup f_m,$$
then we have
$$N_n=M_3\mathop{\cup}_{f}[-1,1]\times F~{\rm and}~\delta_n=(\ga_3\cup\{0\}\times\partial{F})\cap\partial (M_3\mathop{\cup}_{f}[-1,1]\times F).$$

To prove this claim, we focus on the case when $n=1$. The general case follows immediately by induction.

When $n=1$, we have a sutured manifold decomposition
$$(M,\ga)\stackrel{A_1}{\leadsto}(M_{3},\ga_{3})=(M_{1},\ga_{1})\cup (M_{2},\ga_{2}).$$
Write $\partial{A}_1=\al_+\cup\al_-$ so that $\al_{\pm}\subset R_{\pm}(\ga)$. Write $V_{\pm}$ the component of $R_{\pm}(\ga)$ that contains $\al_{\pm}$. We discuss in a few different cases.

When $\al_{+}$ and $\al_{-}$ are non-separating in $V_+$ and $V_-$, respectively, we know that $A_1$ has already been well-groomed, and we just take $S_1=A_1$. Then, we have $(N_1,\delta_1)=(M_3,\ga_3)$. So, for $i=1,...,m$, we simply pick $F_i$ to be an annulus and identify $\theta_i$ with any component of $\partial{F}_i$ but with orientation reversed.

When $\al_+$ is separating, and $\al_-$ is non-separating, we know that $-\al_+$ co-bounds a sub-surface $F_1$ in $V_+$ with part of $\partial V_+$. We can glue $F_1$ to $A_1$ and push it into the interior of $M$. Write the resulting surface $S_1$, then $\partial S_1=\al_-\cup (\partial F_+\backslash\al_+)$, and by Definition \ref{defn: well groomed surface}, $S_1$ is well-groomed. After the decomposition along $A_1$, there is a component of $\ga_3$ corresponding to $\al_+$, which we write $\theta_1$. Then, via $\al_+$, $\theta_1$ is identified with a component of $\partial F_1$ with orientation reversed. It is straightforward to check that
$$N_1=M_3\mathop{\cup}_{\theta_1} [-1,1]\times F_1~{\rm and}~\delta_1=(\ga_3\cup\{0\}\times\partial{F}_1)\cap\partial (M_3\mathop{\cup}_{\theta_1} [-1,1]\times F_1).$$
We can take $F_2,...,F_{m}$ to be annuli just as in the previous case.

When $\al_+$ and $\al_-$ are both separating, assume that $-\al_+$ co-bounds a subsurface $F_+\subset V_+$ together with part of $\partial V_+$, and $-\al_-$ co-bounds a subsurface $F_-\subset V_-$ together with part of $\partial V_-$. We discuss in two cases.

{\bf Case 2.1}. At least one of $F_+$ and $F_-$ has a disconnected boundary. We can glue $F_+$ and $F_-$ to $A_1$ to form the new decomposing surface $S_1$, which is quasi horizontal in the sense of Definition \ref{defn: quasi-horizontal surfaces}. The decomposition of $(M,\ga)$ along $A_1$ yields a sutured manifold $(M_3,\ga_3)$, and two components of $\ga_3$, which we call $\theta_+$ and $\theta_-$, correspond to $\al_+$ and $\al_-$ respectively. As in the above paragraph, $\theta_+$ and $\theta_-$ are identified with a component of $\partial F_+$ and a component of $\partial F_-$, respectively, via $\al_+$ or $\al_-$. Furthermore, if $(N_1,\delta_1)$ is the result of the decomposition of $(M,\gamma)$ along $A_1$, then we know that
$$N_1=M_3\cup [-1,1]\times F_+ \cup [-1,1]\times F_-.$$
Note $(M_3,\gamma_3)$ is taut by construction, and $(N_1,\delta_1)$ is obtained from $(M_3,\gamma_3)$ by gluing product regions $[-1,1]\times F_i$. By Claim 3, which we will prove later, we know that $F_i$ is not a disk, and hence $(N_1,\delta_1)$ is also taut. 

\begin{rem}
Note in this case, the surface $S_1$ might be a horizontal surface instead of just quasi-horizontal (c.f. \cite[Theorem 5.1]{ni2007knot} and \cite[Proposition 6.7]{kronheimer2010knots}). If that happens, then since $(M,\ga)$ is horizontally prime, $S_1$ must be parallel to the boundary. Thus decomposing along $S_1$ peels off a product sutured manifold from $(M,\ga)$. However, the construction of $(N_1,\delta_1)$ is left unchanged, so we don't care if this subtlety happens.
\end{rem}

{\bf Case 2.2}. If both $F_+$ and $F_-$ have connected boundary, then we can simply replace the surface $A_1$ by $-A_1$ to perform the sutured manifold decomposition, and then it falls into Case 2.1. This concludes the proof of Claim 2.

Recall we have
$$(M,\ga)\stackrel{A'}{\leadsto}(M_{3},\ga_{3})=(M_{1},\ga_{1})\cup (M_{2},\ga_{2}).$$
By construction, $(M_1,\gamma_1)$ is reduced, and each component of it is not a product sutured manifold. By \cite[Lemma 2.13]{juhasz2010polytope}, it is also free of essential product disks. Thus, Proposition \ref{prop: decomposition drops the dimension by half} applies and there is a well-groomed surface $S\subset (M_1,\ga_1)$ and a sutured manifold decomposition 
$$(M_1,\ga_1)\stackrel{S}{\leadsto}(M_{4},\ga_{4})$$
so that
$${\rm dim}_{\mathbb{C}}SHI(M_4,\ga_4)\leq\frac{1}{2}{\rm dim}_{\mathbb{C}}SHI(M_1,\ga_1).$$
Note $(M_2,\gamma_2)$ is a product sutured manifold, and $A'$ consists of product annuli, so we know from \cite[Proposition 6.5]{kronheimer2010knots} and Proposition \ref{prop: decompose along product annuli do not increase SHI} that
$${\rm dim}_{\mathbb{C}}SHI(M_4,\ga_4)\leq\frac{1}{2}{\rm dim}_{\mathbb{C}}SHI(M,\ga).$$

From now on, we will call $(M_4,\ga_4)\sqcup (M_2,\ga_2)$ still $(M_4,\ga_4)$, so we have a well-groomed surface $S\subset (M_3,\ga_3)$ and a taut sutured manifold decomposition
$$(M_3,\ga_3)\stackrel{S}{\leadsto}(M_{4},\ga_{4})$$
with
$${\rm dim}_{\mathbb{C}}SHI(M_4,\ga_4)\leq\frac{1}{2}{\rm dim}_{\mathbb{C}}SHI(M,\ga).$$

Next, we need to modify $(N_n,\delta_n)$ for the purpose of further discussions. Suppose for some $i\in\{1,...,m\}$, the surface $F_i$ has a connected boundary. Then, $\partial{F}_i$ is necessarily identified with $\theta_i$, with orientation reversed.

{\bf Claim 3}. In the proof of Claim 2, if the surface $F_i$ has a connected boundary, then it is not a disk.

Suppose, on the contrary, that $F_i$ is a disk. Then recall, by construction, the surface $F_i$ is a subsurface of $R(\ga)\subset \partial M$, and there is a component $A_j$ of $A'$ and a boundary component $\al$ of $A_j$ so that $\partial{F_i}=-\al$. Let $\al'$ be the other boundary component of $A_j$, then we know that $A_j\cup F_i$ is a disk whose boundary is $\al'$. Since $(M,\ga)$ is taut, we know that $\al'$ also bounds a disk $D\subset R(\ga)$. Then, $A_j\cup F_i\cup D$ is a $2$-sphere. Since $(M,\ga)$ is taut, this $2$-sphere bounds a $3$-ball, and hence $A_j$ is a trivial product annulus, which contradicts the way we choose $A_j$.

Now we explain how to modify $(N_{n},\delta_n)$: Suppose for some $i\in\{1,...,m\}$, $F_i$ has a connected boundary and is glued to a component of $\ga_1$. By Claim 3, $F_i$ is not a disk, so we can pick a non-separating curve $\be_i$ on $F_i$. Then, 
$$[-1,1]\times\be_i\subset [-1,1]\times F_i\subset (N_n,\delta_n)$$
is a product annulus, which is also a well-groomed surface. Let $A''$ be the union of all such product annuli, one for each $F_i$ so that it has a connected boundary that is glued to a component $\ga_1$. Since $(M_1,\ga_1)$ is reduced, we can argue in the same way as in the proof of Claim 1 and conclude the following.

{\bf Claim 4}. We have
$$|A''|\leq|\ga_1|< 6\cdot 2^{k+1}.$$

We have a sutured manifold decomposition
$$(N_n,\delta_n)\stackrel{A''}{\leadsto}(N_n',\delta_n').$$ By construction, we know that there are connected compact oriented surfaces $F_i'$, which are either $F_i$ or $F_i$ cut open along $\be_i$, so that the following is true. (Recall $\theta_i$ are the components of $\ga_3$.)

(a') Each $F_i'$ has at least two boundary components.

(b') For $i=1,...,m$, there is an orientation reversing embedding
$$f'_i:\theta_i\hookrightarrow \partial{F}'_i.$$

(c') Write 
$$F'=F'_1\cup...\cup F'_m ~{\rm and}~ f'=f'_1\cup...\cup f'_m,$$
then we have
$$N'_n=M_3\mathop{\cup}_{f'}[-1,1]\times F'~{\rm and}~\delta'_n=(\ga_3\cup\{0\}\times\partial{F'})\cap\partial (M_3\mathop{\cup}_{f'}[-1,1]\times F').$$

Next, we want to extend the well-groomed surface $S\subset (M_3,\ga_3)$ to a well-groomed surface $S'$ on $(N_n',\delta_n')$. To do this, we extend $S$ across all $[-1,1]\times F_i'$ as follows: For $i\in\{1,...,m\}$, let $A(\theta_i)$ be the annular neighborhood of $\theta_i\subset \ga_3$. If $S\cap A(\theta_i)=\emptyset$, then we are already done. If $S\cap A(\theta_i)$ consists of parallel copies of $\theta_i$, then we simply glue the same number of copies of $F_i$ to $S$. If $S\cap A(\theta_i)$ is neither empty nor a collection of simple closed curves, then $S$ being well-groomed implies that $S\cap \theta_i$ is a finite set of points of the same signs. Write $\sigma_i'$ the component of $\partial{F}_i'$ so that $\sigma_i'=f'(\theta_i)$. Let $\tau_i'\subset F_i'$ be a disjoint union of parallel arcs so that each component of $\tau_i'$ has one end point on $\sigma_i'$ and the other end point on $\partial F_i'\backslash \sigma_i'$, which is non-empty by condition (a'). We further require that
$$\tau_i'\cap \sigma_i'=f'(S\cap \theta_i).$$
Then, we can glue $[-1,1]\times \tau_i'$ to $S$ along $[-1,1]\times (\tau_i'\cap \sigma_i')$. 

Performing this operation for all $i$, we obtain a surface $S'\subset (N_n',\delta_n')$.

To show that $S'$ is well-groomed, first it is straightforward to check that $\partial S'$ is essential in $H_1(\partial N_n')$, and condition (1) in Definition \ref{defn: well groomed surface} is satisfied by the construction of $S'$. To show that condition (2) also holds, suppose $V'_n$ is a component of $R(\delta_n')$. By construction, there exists a component $V_3$ of $R(\ga_3)$, and surfaces $F_{i_1}$,..., $F_{i_{l}}$ so that
$$V'_n=V_3\cup F_{i_1}\cup...\cup F_{i_l}.$$
It is crucial that our construction of $(N_n',\delta_n')$ makes sure that one component of $R(\delta_n')$ contains only one component of $R(\ga_3)$. If $\partial S\cap V_3$ consists of parallel and coherently oriented non-separating simple closed curves, then
$$\partial S'\cap V_n'=\partial S\cap V_3.$$
If $\partial S\cap V_3$ consists of parallel and coherently oriented properly embedded simple arcs, then each arc intersects two components of $\partial V_3$, say $\theta_{i_1}$ and $\theta_{i_2}$. Then we know that each component of $\partial S'\cap V_n'$ is an arc by gluing three pieces, $\tau_{i_1}'$, $\tau_{i_2}'$, and a component of $\partial S\cap V_3$, together. Hence by Definition \ref{defn: well groomed surface}, $S'$ is indeed well-groomed.

Now there are two decompositions:
$$(M_3,\ga_3)\stackrel{S}{\leadsto}(M_{4},\ga_{4}){\rm~and~}(N_n',\delta_n')\stackrel{S'}{\leadsto}(N_{n+1},\delta_{n+1}).$$

{\bf Claim 5}. We have
$$SHI(M_4,\ga_4)\cong SHI(N_{n+1},\delta_{n+1}).$$

To prove this claim, for any $i\in\{1,...,m\}$, let $A(\theta_i)$ be the annular neighborhood of $\theta_i\subset \ga_3$. If $S\cap A(\theta_i)=\emptyset$, then $A(\theta_i)$ survives in $M_4$, and $\theta_i$ is a component of $\ga_4$. To obtain $N_{n+1}$ from $M_4$, we need to glue $[-1,1]\times F_i'$ to $M_4$ along $\theta_i$. If $S\cap A(\theta_i)$ consists of parallel copies of $\theta_i$, then $(N_{n+1},\delta_{n+1})$ is obtained from $(M_4,\ga_4)$ by gluing a few copies of $[-1,1]\times F_i'$. By \cite[Proposition 6.7]{kronheimer2010knots}, gluing product regions will not change the sutured instanton Floer homology. Hence, it remains to deal with the case when $S\cap \theta_i\neq\emptyset$. Suppose
$$S\cap \theta_i=\{p_1,...,p_{s_i}\},$$
where $p_1,...,p_{s_i}$ are labeled according to the orientation of $\theta_i$. Let $\theta_i''$ be the part of $\theta_i$ from $p_{s_i}$ to $p_1$. Then, $\theta_i''$ does not contain any other $p_j$. Recall there is a collection of arcs $\tau_i'\subset F_i'$. Note $F_i'\backslash \tau_i'$ consists of a few disks and a large piece $F_{i}''$ that contains most of $F_i'$. Furthermore, there is a component $\sigma_i''$ of $\partial F_i'\backslash \tau_i'\subset \partial F_{i}''$ so that $\sigma_i''=f'(\theta_i'')\subset \partial F_i''$. It is straightforward to check that, to obtain $N_{n+1}$ from $M_4$, we need to glue $[-1,1]\times F_i''$ to $M_4$ along $[-1,1]\times \sigma_i''$. Note $\sigma_i''$ is an arc, so topologically $[-1,1]\times \sigma_i''$ is a disk which intersects the suture $\ga_4$ along an arc $\theta_i''$, and intersects the suture $\{0\}\times\partial{F}_i''$ of the product sutured manifold $([-1,1]\times F_{i}'',\{0\}\times \partial F_i'')$ along another arc $\{0\}\times\sigma_i''$. Hence, this gluing coincides with the setting of attaching a contact $1$-handle to the disjoint union
$$(M_4,\ga_4)\sqcup([-1,1]\times F_{i}'',\{0\}\times \partial F_i''),$$
in the sense of Baldwin and Sivek \cite{baldwin2016instanton}. Since both disjoint union with a product sutured manifold and attaching a contact $1$-handle do not change sutured instanton Floer homology, we conclude that
$$SHI(M_4,\ga_4)\cong SHI(N_{n+1},\delta_{n+1}).$$

Finally, we are ready to finish the induction. We have a sequence of decompositions:
\begin{equation}\label{eq: sequence of decompositions}
(M,\ga)\stackrel{S_1}{\leadsto}(N_{1},\delta_{1})\stackrel{S_2}{\leadsto}...\stackrel{S_n}{\leadsto}(N_{n},\delta_{n})\stackrel{A''}{\leadsto}(N_n',\delta_n')\stackrel{S'}{\leadsto}(N_{n+1},\delta_{n+1}).
\end{equation}
We know from Claim 5 that
\beq{\rm dim}_{\mathbb{C}}SHI(N_{n+1},\delta_{n+1})&={\rm dim}_{\mathbb{C}}SHI(M_4,\ga_4)\\
&\leq \frac{1}{2}{\rm dim}_{\mathbb{C}}SHI(M_3,\ga_3)\\
&\leq \frac{1}{2}{\rm dim}_{\mathbb{C}}SHI(M,\ga)\\
&<2^{k}.
\eeq
Thus, the inductive hypothesis applies on $(N_{n+1},\delta_{n+1})$, and there is a taut foliation $\mathfrak{F}'$ of depth at most $2^{k+5}$ on $(N_{n+1},\delta_{n+1})$. We now go through decomposition (\ref{eq: sequence of decompositions}) to construct a taut foliation on $(M,\ga)$. First, each $S_i$ is either well groomed or quasi horizontal by Claim 2, and $n<6\cdot 2^{k+1}$ by Claim 1. Second, each component of $A''$ is well-groomed and, since the components of $A''$ are contained in sufficiently disjoint regions of $N_n$, decomposing along some subset of $A''$ will keep each of the rest components of $A''$ being well-groomed. By Claim 4, $|A''|< 6\cdot 2^{k+1}$. Finally, there is one last decomposition along the well-groomed surface $S'$, so by Theorem \ref{thm: well groomed decomposition increase depth by 1}  and Lemma \ref{lem: glue taut foliation through quasi-horizontal surfaces}, there is a taut foliation $\mathfrak{F}$ on $(M,\ga)$ of depth at most
$$6\cdot 2^{k+1}+6\cdot 2^{k+1}+1+2^{k+5}<2^{k+6}.$$
Hence, the inductive step is completed, and we finish the proof of Theorem \ref{thm: small depth taut foliations}.
\epf

\bcor\label{cor: SFH}
Suppose $(M,\ga)$ is a taut balanced sutured manifold with $H_2(M)=0$, and
$${\rm rk}_{\intg_2}(SFH(M,\ga))<2^{k+1},$$
then $(M,\ga)$ admits a taut foliation of depth at most $2^{k+6}$.
\ecor

\bpf
The proof of Theorem \ref{thm: small depth taut foliations} applies verbatim.
\epf

The above corollary can be used to prove the following, which gives a partial answer to \cite[Question 9.14]{juhasz2008floer}.

\bcor\label{cor: HFK}
Let $K$ be a knot in a rational homology $3$-sphere $Y$ so that the knot complement $Y(K)$ is irreducible. Suppose further
that $k$ is a positive integer so that
$${\rm rk}_{\intg_2}\widehat{HFK}(Y,K,g(K))<2^{k}.$$
Then, $Y\backslash N(K)$ admits a taut foliation of depth at most $2^{k+5}$ transverse to the boundary of $N(K)$.
\ecor
\bpf
We have a sutured manifold $Y\backslash N(K)$, with toroidal suture. Take $S$ to be a minimal genus rational Seifert surface of $K$, then we have a sutured manifold decomposition
\begin{equation}\label{eq: decomposing knot complement in rational homology sphere}
(Y\backslash N(K))\stackrel{S}{\leadsto}(M,\ga),
\end{equation}
and we know that
$$SFH(M,\ga)\cong \widehat{HFK}(Y,K,g(K)).$$
Thus, Corollary \ref{cor: SFH} applies and we obtain a taut foliation on $(M,\ga)$ of depth at most $2^{k+5}$. We can further glue it along the decomposition (\ref{eq: decomposing knot complement in rational homology sphere}), which will not increase the depth of the taut foliation. Hence, we are done.
\epf

\bpf[Proof of Corollary \ref{cor: small depth taut foliations}]
On the knot complement $S^3(K)=S^3\backslash N(K)$, we can pick $\Ga_{\mu}$ to be a suture consisting of two meridians. From Corollary 4.2 in Kronheimer and Mrowka \cite{kronheimer2010instanton}, we know that
$${\rm dim}_{\mathbb{C}}SHI(S^3(K),\Ga_{\mu})< 1+2^{k+2}.$$
Pick a minimal genus Seifert surface $S$ of $K$, we know that there is a decomposition
$$(S^3(K),\Ga_{\mu})\stackrel{S}{\leadsto}(M,\ga),$$
and by \cite[Lemma 5.7 and Proposition 5.11]{kronheimer2010knots}, we know that
$${\rm dim}_{\mathbb{C}}SHI(M,\ga)<2^{k+1}.$$
Hence, we can apply Theorem \ref{thm: small depth taut foliations}, and there is a taut foliation on $(M,\ga)$ of depth at most $2^{k+6}$. Note we can also regard $S^3(K)$ as a sutured manifold with toroidal sutures, and decomposing $S^3(K)$ along $S$ also gives rise to $(M,\ga)$. So, we can glue the just obtained taut foliation on $(M,\ga)$ along this later decomposition to conclude the proof of Corollary \ref{cor: small depth taut foliations}. 
\epf


\bibliography{Index}
\end{document}